\documentclass[11pt]{article}
\usepackage{amsfonts}
\usepackage{latexsym,amsmath}
\usepackage{amssymb,array}
\RequirePackage[dvips]{graphicx}
\usepackage{graphics}
\usepackage{tikz}
\usepackage{adjustbox}
\usepackage{caption}
\usepackage{calc}
\usepackage{graphics}
\usepackage{longtable}
\usepackage{pdflscape}
\usepackage{xcolor}

\usepackage{authblk}
\usepackage[colorlinks=true]{hyperref}
\hypersetup{urlcolor=blue, citecolor=red}
\makeatletter \oddsidemargin 0in \evensidemargin 0in \textwidth
16cm\textheight 20cm 
\begin{document}
\newcommand{\la}{\lambda}
\newcommand{\eq}{\Leftrightarrow}
\newcommand{\mf}{\mathbf}
\newcommand{\ri}{\Rightarrow}
\newtheorem{t1}{Theorem}[section]
\newtheorem{d1}{Definition}[section]
\newtheorem{n1}{Notation}[section]
\newtheorem{c1}{Corollary}[section]
\newtheorem{l1}{Lemma}[section]
\newtheorem{r1}{Remark}[section]
\newtheorem{e1}{Counterexample}[section]
\newtheorem{p1}{Proposition}[section]
\newtheorem{cn1}{Conclusion}[section]
\newtheorem{il1}{Illustration}[section]
\renewcommand{\theequation}{\thesection.\arabic{equation}}
\pagenumbering{arabic}
\title {Stochastic Comparisons of Random Extremes from Non-Identical Random Variables}

\author[1]{Amarjit Kundu}
\author[2]{Shovan Chowdhury\footnote{Corresponding author e-mail:meetshovan@gmail.com; shovanc@iimk.ac.in}}
\author[1]{Bidhan Modok}

\affil[1] {Department of Mathematics, Raiganj University, West Bengal, India}
\affil[2] {Quantitative Methods and Operations Management Area, Indian Institute of Management Kozhikode, India}

\maketitle
\begin{abstract}
We propose some new results on the comparison of the minimum or maximum order statistic from a random number of non-identical random variables. Under the non-identical set-up, with certain conditions, we prove that random minimum (maximum) of one system dominates the other in hazard rate (reversed hazard rate) order.
Further, we prove variation diminishing property (\cite{ka}) for all possible restrictions to derive the new results. 

\end{abstract}
%----------------------------------------------------Introduction-------------------------
\setcounter{section}{0}
\section{Introduction}
\setcounter{equation}{0}
\hspace*{0.3in} The problem of stochastic comparison of extremes (maximum or minimum) has always been an important topic in a variety of disciplines when sample size is considered fixed or random. In a fixed sample size scenario, many researchers have compared series (minimum order statistic) or parallel (maximum order statistic) system when the components are identical or non-identical and are ordered in some fashion (see~\cite{tr11},\cite{kun2},~\cite{bar3},~\cite{kun5},~\cite{ch}). Moreover, readers may refer to~\cite{bz11.1} and the references therein for a comprehensive review til 2013. However, in the field of finance and risk management, reliability optimization and life testing experiment, actuarial science, biological Science, engineering and manufacturing, hydrology and others, sample size may be unobserved, or may depend on the occurrence of unsystematic events which makes it random. Following examples justify the use of random sample size in practice.   \\
\hspace*{0.3in} Lifetime experiments are always censored and hence sample size can hardly be fixed. In actuarial science, the performance of a portfolio of insurance contracts requires stochastic modeling of the largest claims while the number of claims is random. In biological and agricultural experiments, some elements of the sample under consideration may die for unknown reasons making sample size random. In finance and risk management, one may be interested in the minimum and maximum loss from a portfolio of loans or securities with random number of assets or liabilities. The use of random sample size arises naturally in reliability engineering to model lifetime of a parallel or series system. One can think of a situation where failure (of a device for example) occurs due to the presence of an unknown number (random) of initial defects of the same kind (a number of semiconductors from a defective lot, for example) where each defect can be detected only after causing failure, in which case it is repaired perfectly. In reliability optimization, one system can dominate the other under certain conditions when the components are chosen randomly from two batches (\cite{di},~\cite{ha},~\cite{na},~\cite{nav}). In hydrological applications (\cite{ku}) concerning the distribution of annual maximum rainfall or flood, it may be assumed that the number of storms or floods per year is not constant; in traffic engineering applications (\cite{sr}) velocity distribution of vehicles on a particular spot on a road or a highway is obtained after recording velocities of a random number of vehicles pass through the point; in transportation theory (\cite{sw}) modeling accident-free distance of a shipment of explosives require random number of explosives which are defective and lead to explosion and cause an accident; in biostatistical applications (\cite{kou} concerning the distribution of minimum time for clonogenic (carcinogenic) cell to produce a detectable cancer mass, number of clonogens (carcinogenic cells) left active after a cancer treatment is always random.\\ 
\hspace*{0.3in} Situations where sample size is random, one can find results concerning stochastic comparison of minima and maxima for independent and identically distributed ($iid$) random variables. For $iid$ random variables $X_i,~i=1,2,3...$ having common distribution function $F$ and for a discrete random variable $N$,~\cite{sh} showed that there exist bounds, both lower and upper, of the failure rate functions of random minimum $X_{1:N}=min(X_1,X_2,…,X_N)$ and random maximum $X_{N:N}=max(X_1,X_2,…,X_N)$, in terms of failure rate function of $F$. In a similar set up,~\cite{ba} showed that $X_{1:N}~(X_{N:N})$ is smaller than $Y_{1:N}~(Y_{N:N})$~in terms of convex transform order, star order and super additive order if $X_1$ is smaller than $Y_1$ with respect to the same ordering. Likelihood ratio ordering was also established under the same set up;~\cite{li} proved similar results for right spread ordering and increasing convex ordering between $X_{N:N}$ and $Y_{N:N}$ and for total time on test transform ordering and increasing concave ordering between $X_{1:N}$ and $Y_{1:N}$.~\cite{ah} proved the reversed result which states that if there exists right spread ordering and increasing convex ordering between $X_{1:N}$ and $Y_{1:N}$ and total time on test transform ordering and increasing concave ordering between $X_{N:N}$ and $Y_{N:N}$, then the same ordering exists between $X_1$ and $Y_1$. For $iid$ random variables and for discrete random variables $N_1$ and $N_2$,~\cite{sw} showed that if $N_1\leq_{rh}N_2,$ then $X_{1:N_1}\geq_{hr}X_{1:N_2}$ and $X_{N_1:N_1}\leq_{rh}X_{N_2:N_2}$. It was also shown that if $N_1\leq_{hr}N_2,$ then $X_{1:N_1}\geq_{rh}X_{1:N_2}$ and $X_{N_1:N_1}\leq_{hr}X_{N_2:N_2}$. The authors also proved that if $N_1 \leq_{lr}N_2,$ then $X_{1:N_1}\geq_{lr}X_{1:N_2}$ and $X_{N_1:N_1}\leq_{lr}X_{N_2:N_2}$. For an overview of stochastic comparisons with random sample size, see~\cite{ns}.
\\ So, the natural question arises whether random extremes can be compared for independent and non-identically distributed ($inid$) random variables. To the best of our knowledge, there has been no attention paid to the case of non-identical set-up which motivates us to compare random extremes of $inid$ random variables with suitable stochastic ordering principle. To this endeavor, we have proven variation diminishing property (\cite{ka}) for all possible restrictions which fills some theoretical gap in the literature. \\
\hspace*{0.3in} One may think of applying the key findings of the paper in the following real life situations. \\ Fiber composites of many products are made of carbon fiber, breaking strength of which is taken to be one of the significant quality characteristics in engineering design and material science. Materials of high strength are always desirable for human safety and well-being. Hence, it is essential to monitor the quality of materials produced in any manufacturing process. Suppose there are two independent quality control stations in a production line inspecting random number of products for defects. It is worthwhile to compare the random maximum of the breaking strengths of carbon fiber taken from two independent lots. Assume that the maximum breaking strength from the first lot of $n$ carbon fibers is more than the same from the second in reversed hazard rate ordering. Theorem~\ref{th3} guarantees that the same result will hold under certain conditions for random number of items. \\ Let us consider two independent queueing systems with one server each. One may wish to compare the service time of customers stochastically in order to measure the operational efficiency of the servers. Time taken to serve the $i^{th}$ customer ($i=1,2,...,n$) by the $j^{th}$ server ($j=1,2$) can be assumed to be $inid$ random variable. It is known that the minimum of the service time of $n$ customers served by server 1 is less than that of the same by server 2 in hazard rate ordering. In practice, queue size is random and hence it is important to answer what happens to the comparison when the number of customers ($N$) served by each server is random and follows any discrete distribution $p(n)$ with common parameter(s) with $\mathcal{N^+}$ as support. According to Theorem~\ref{th1}, minimum of the service time of $N$ customers served by server 1 is also less than that of the same by server 2 in hazard rate ordering. The result confirms that the same ordering relationship will hold for random number of customers under certain condition.  \\ Suppose two portfolios of insurance contracts need to be compared stochastically in terms of the maximum claim amounts (loss to the firm). The actual number of claims is supposed to be always random and the portfolio managers are interested to know which portfolio incurs maximum loss. If it is known that the maximum loss from the first portfolio of $n$ claims is more than the same from the second in reversed hazard rate ordering, according to Theorem~\ref{th3}, the same result will hold under certain conditions for random claim size $N$ defined with support $\mathcal{N^+}.$ \\ 
\hspace*{0.3 in} The organization of the paper is as follows. Relevant definitions and existing results are given in Section 2. In section 3, new lemmas are derived and the results related to hazard rate order and reversed hazard rate order between two random extremes of $inid$ random variables are established with examples. Section 4 concludes the paper. 
%%%%%%%%%%%%%%%%%%%%%%%%%%%%%%%%%%%%%%%%%%%%%%%%%%%%%%%%%%%%%%%%%%%%%%%%%%%%%%%%%%%%%%%%%%%%%%%%%%%%%%%%%%%%%%%%%%%%%%%%%%%%%%%%%%%%%%%%%%%%
\\\hspace*{0.3 in}Throughout the paper, the word increasing (resp. decreasing) and nondecreasing (resp. nonincreasing) are used interchangeably. The notations $\mathcal{R}_+$ and $\mathcal{N}^+$ denote the set of positive real numbers $\{x:0<x<\infty\}$ and natural numbers $\{x:x=1,2,...,\infty\}$. The random variables considered in this paper are all nonnegative.

%\\\hspace*{0.3in} 
\section{Preliminaries}
\setcounter{equation}{0}
\hspace*{0.3 in} For an absolutely continuous random variable $X$, we denote the cumulative distribution function by $F_X(\cdot)$, the hazard rate function by $r_X(\cdot)$ and the reversed hazard rate function $\overline r_X(\cdot)$. The survival or reliability function of the random variable $X$ is written as $\bar F_X(\cdot)=1-F_X(\cdot)$. 
\\\hspace*{0.3 in} In order to compare different order statistics, stochastic orders are used for fair and reasonable comparison. In literature different kinds of stochastic orders are developed and studied. The following well known definitions can be obtained in~\cite{shak}.
\begin{d1}\label{de1}
Let $X$ and $Y$ be two absolutely continuous random variables with respective supports $(l_X,u_X)$ and $(l_Y,u_Y)$, where $u_X$ and $u_Y$ may be positive infinity, and $l_X$ and $l_Y$ may be negative infinity. Then, $X$ is said to be smaller than $Y$ in
\begin{enumerate}
%\item likelihood ratio (lr) order, denoted as $X\leq_{lr}Y$, if 
%$$\frac{f_Y(t)}{f_X(t)}\;\text{is increasing in} \,t\in(l_X,u_X)\cup(l_Y,u_Y);$$
\item likelihood ration order, denoted as $X\leq_{lr}Y$; if $$\frac{f_Y(t)}{f_X(t)}\;\text{is increasing in}\, t \in (l_x,u_x)\cup(l_y,u_y),$$
\item hazard rate (hr) order, denoted as $X\leq_{hr}Y$, if $$\frac{\bar F_Y(t)}{\bar F_X(t)}\;\text{is increasing in}\, t \in (-\infty,max(u_X,u_Y)),$$
 which can equivalently be written as $r_X(t)\geq r_Y(t)$ for all $t$;
 \item reversed hazard rate (rh) order, denoted as $X\leq_{rh}Y$, if $$ \frac{F_Y(t)}{ F_X(t)}\;\text{is increasing in}\, t \in(min(l_X,l_Y),\infty),$$
 which can equivalently be written as $\overline r_X(t)\leq \overline r_Y(t)$ for all $t$;
 \item usual stochastic (st) order, denoted as $X\leq_{st}Y$, if $\bar F_X(t)\leq \bar F_Y(t)$ for all \\$t\in (-\infty,\infty).$
$\hfill\Box$ 
 %\item up shifted reversed hazard rate (rhr $\uparrow$) order, denoted as $X\leq_{rhr \uparrow}Y$, if $X-x\leq_{rhr}Y,$ for all $x\geq 0$, or equivalently, if 
 %$$\frac{F_Y(t)}{F_X(t+x)} \;\text{is increasing in}\; t\in(l_X,\infty),\;\text{for all} \;x\geq 0.$$ $\hfill\Box$
\end{enumerate}
\end{d1}
\hspace*{0.3in} The next two theorems due to~\cite{shak} are used to prove the key findings of the paper.
\begin{t1}\label{theo11}
 $\left(\text{Theorem 1.B.28, Page 31}\right)\;$
If $X_1,\ X_2\ldots,\ X_m$ are independent random variables, then $X_{(k:m-1)}\geq_{hr}X_{(k:m)}$ for $k=1,2,\ldots, m-1.$
\end{t1}
\begin{t1}\label{theo12}
 $\left(\text{Theorem 1.B.57, Page 41}\right)\;$
Let $X_1,\ X_2\ldots,\ X_m$ be independent random variables. If $X_m\leq_{rh}X_j$ for all $j=1,2,\ldots, m-1$, then $X_{(k-1:m-1)}\leq_{rh}X_{(k:m)}$ for $k=2,3,\ldots, m$.
\end{t1}
\hspace*{0.3in} The notion of Totally Positive of order $2$ ($TP_2$) and Reverse Regular of order $2$ ($RR_2$) are of great importance in various fields of applied probability and statistics (see~\cite{ka},~\cite{ka1,ka2}). It can be recalled that a non-negative function $f:\Re^2\longmapsto\Re_+$ is said to be $TP_2$ if
\begin{equation}\label{eq1a}
	f(u_1,v_1)f(u_2,v_2)-f(u_2,v_1)f(u_1,v_2)\geq 0,
\end{equation}
for all $u_1\leq u_2$ and $v_1\leq v_2$. This $f$ is said to be $RR_2$ if the inequality in (\ref{eq1a}) is reversed. This can equivalently be written as follows\\
$f$ is $TP_2$ ($RR_2$) if $f(u,v_2)/f(u,v_1)$ is increasing (decreasing) in $u$ for all $v_1\leq v_2$. \\
\hspace*{0.3 in} The following Theorem on variation diminishing property is proved by~(\cite{ka}) for some definite restrictions on $f(x,y)$ and $K(x,y)$.   
\begin{t1}\label{theo13}
 $\left(\text{Theorem 11.2, Page 324-25}\right)\;$~Let $K\left(x,y\right)>0$, defined on $X \times Y,$ be $TP_2$, where $X$ and $Y$ are the subsets of the real line. Assume that a function $f\left(x,y\right)$ is such that
\begin{enumerate}
\item[(i)] for each $x$, $f\left(x,y\right)$ changes sign at most once, and from negative to positive values, as $y$ traverses $Y;$
\item[(ii)] for each $y$, $f\left(x,y\right)$ is increasing in $x$;
\item[(iii)] $g\left(x\right)=\int_{Y} f\left(x,y\right)K\left(x,y\right)d\mu(y)$ exists absolutely and defines a continuous function on $x$.
\end{enumerate}
Then $g\left(x\right)$ changes sign at most once, and from negative to positive values.$\hfill\Box$
\end{t1}
The Propositions in the next section are based on Theorem~\ref{theo13} yielding results for all possible restrictions on $f(x,y)$ and $K(x,y)$ where $y$ is defined on $\mathcal{N^+}.$ The propositions contribute to the literature in the field of inequalities and are the building blocks for the key results of the paper. Some counter examples are also provided in the section to show that for other combinations of conditions as given in the propositions, the results do not hold.
%In the following diagrams we present a chain of implications of the stochastic orders (cf. Shaked and Shanthikumar \cite{shak1}, and Lillo \emph{et al.} \cite{lns1}):
%\vspace{0.17 in}
%\[ \xrightarrow{\hspace*{3cm}} \]
%\\\hspace*{1.7 in}$~~~~~~ X\leq_{hr}Y$
%\\\hspace*{1.7 in}$~~\nearrow ~~~~~~~~~~~~\searrow$
%\\\hspace{6 in} $~~~~~~~~~~~~~~~~~~~~~~~X\leq_{lr}Y~\rightarrow~~~~~~~~~~~~~~X\leq_{st}Y$

%\hspace*{1.7 in}$~~\searrow ~~~~~~~~~~~~\nearrow$
%\\\hspace{6 in}$~~~~~~~~~~~~~~~~~~~X\leq_{rhr\uparrow}Y~\rightarrow~~ X\leq_{rhr}Y$\\

\section{Main Results}
This section is divided into two subsections. In the first subsection four propositions, mentioned earlier, has been proved. One counterexample is also provided to show that, for other combinations of conditions, mentioned in the propositions, no conclusion can be drawn. In the second subsection results related to comparisons of two maximum and minimum order statistics for random number of components, are provided.
\subsection{Some Propositions}
In this subsection readers will find the proof of Proposition~\ref{le1}. The other three propositions can be proved in a similar line and hence are omitted.
\begin{p1}\label{le1} For any positive integer $n$ and positive real number $x$, let $K_n(x)>0$ be a $RR_2$ function in $n\in \mathcal{N}^+$ and $x\in \mathcal{R}_+$. Assume that any function $f_n(x)$ be such that 
\begin{itemize}
\item [i)] for each $x\in \mathcal{R}_+$, $f_n(x)$ changes sign at most once and if the change occurs, it is from negative to positive as $n$ traverses in $\mathcal{N}^+$;
\item [ii)] for each $n\in \mathcal{N}^+$, $f_n(x)$ decreases in $x$;
\item [iii)] $w(x)=\sum_{n=1}^{\infty}{f_n(x)K_n(x)}$ converges absolutely and defines a continuous function of $x$.
\end{itemize}
Then $w(x)$ changes sign at most once and if the change occurs, it is from positive to negative.
\end{p1} 
\textbf{Proof:}  To prove the result, it is sufficient to show that for any $x_0\in\mathcal{R}_+$, $w(x_0)=0$ implies $w(x)\geq 0$ for all $x \leq x_0.$\\
For any $x_0\in\mathcal{R}_+$ let us assume that $w(x_0)=0$. Then there must exist one positive integer $n_0\in \mathcal{N^+}$ such that $f_n(x_0)\geq 0$ for all $n\geq n_0$ and $f_n(x_0)\leq 0$ for all $n\le n_0.$ For any $x\in \mathcal{R}_+$ let us consider
\begin{eqnarray}\label{eq1}
\frac{w(x)}{K_{n_0}(x)} &=& \frac{w(x)}{K_{n_0}(x)}- \frac{w(x_0)}{K_{n_0}(x_0)} \quad \left[\text{since} \; w(x_0)=0\right]\nonumber\\
												&=& \sum_{n=1}^{\infty}{\frac{K_n(x)f_n(x)}{K_{n_0}(x)}-\sum_{n=1}^{\infty}{\frac{K_n(x_0)f_n(x_0)}{K_{n_0}(x_0)}}}\nonumber\\
												&=& \sum_{n=1}^{\infty}{\left[\frac{K_n(x)f_n(x)}{K_{n_0}(x)}-\frac{K_n(x)f_n(x_0)}{K_{n_0}(x)}+\frac{K_n(x)f_n(x_0)}{K_{n_0}(x)}-\frac{K_n(x_0)f_n(x_0)}{K_{n_0}(x_0)}\right]}\nonumber\\
												&=& \sum_{n=1}^{\infty}{f_n(x_0)\left[\frac{K_n(x)}{K_{n_0}(x)}-\frac{K_n(x_0)}{K_{n_0}(x_0)}\right]}+\sum_{n=1}^{\infty}{\left[f_n(x)-f_n(x_0)\right]\frac{K_n(x)}{K_{n_0}(x)}}.
\end{eqnarray}
Now two cases may arise.\\
\underline{Case-1}: Let $x<x_0$ and $n>n_0$. Now as $K_n(x)$ is $RR_2$ in $n$ and $x$, it can be written that $K_n(x_0)K_{n_0}(x)\leq K_{n_0}(x_0)K_n(x),$ for all $x<x_0$ and $n>n_0,$ which gives 
$$\frac{K_n(x)}{K_{n_0}(x)}-\frac{K_n(x_0)}{K_{n_0}(x_0)}\geq 0.$$ 
Since $n>n_0$ implies $f_n(x_0)>0$, the first term of (\ref{eq1}) is positive. Now, for each $n\in\mathcal{N}^+,$ $f_n(x)$ is decreasing in $x$ and $x<x_0$, yielding $f_n(x)-f_n(x_0)\geq 0$. Hence, from the fact that, for all $n\in \mathcal{N}^+$ and $x\in \mathcal{R}_+$ $K_n(x)>0$, it can be concluded that the second term of (\ref{eq1}) is also positive.\\
\underline{Case-2}: Now, let $x<x_0$ and $n<n_0$. Since $K_n(x)$ is $RR_2$ in $n$ and $x$, then $K_n(x)K_{n_0}(x_0)\leq K_{n_0}(x)K_n(x_0)$, which in turn gives 
$$\frac{K_n(x)}{K_{n_0}(x)}-\frac{K_n(x_0)}{K_{n_0}(x_0)}\leq 0.$$
Since $n<n_0$ implies $f_n(x_0)<0,$ the first term of ($\ref{eq1}$) is positive. Again, for each $n,$ $f_n(x)$ is decreasing in $x$, yielding $f_n(x)-f_n(x_0)\geq 0,$for all $x<x_0.$ Therefore, from the fact that for all $n, x,$ $K_n(x)>0$, it can be concluded that the second term of (\ref{eq1}) is also positive and thus $\frac{w(x)}{K_{n_0}(x)}$ is positive.\\
Hence, for $x<x_0,$ $w(x)>0$ in both the cases. Thus the result is proved.$\hfill\Box$\\

\begin{p1} \label{le2}
For any positive integer $n$ and positive real number $x$, let $K_n(x)>0$ be a $RR_2$ function in $n\in \mathcal{N}^+$ and $x\in \mathcal{R}_+$. Assume that any function $f_n(x)$ be such that 
\begin{itemize}
\item [i)] for each $x\in \mathcal{R}_+$, $f_n(x)$ changes sign at most once and if the change occurs, it is from positive to negative as $n$ traverses in $\mathcal{N}^+$;
\item [ii)] for each $n\in \mathcal{N}^+$, $f_n(x)$ increases in $x$;
\item [iii)] $w(x)=\sum_{n=1}^{\infty}{f_n(x)K_n(x)}$ converges absolutely and defines a continuous function of $x$.
\end{itemize}
Then $w(x)$ changes sign at most once and if the change occurs, it is from negative to positive.
\end{p1} 

\begin{p1}\label{le3}
For any positive integer $n$ and positive real number $x$, let $K_n(x)>0$ be a $TP_2$ function in $n\in \mathcal{N}^+$ and $x\in \mathcal{R}_+$. Assume that any function $f_n(x)$ be such that 
\begin{itemize}
\item [i)] for each $x\in \mathcal{R}_+$, $f_n(x)$ changes sign at most once and if the change occurs, it is from positive to negative as $n$ traverses in $\mathcal{N}^+$;
\item [ii)] for each $n\in \mathcal{N}^+$, $f_n(x)$ decreases in $x$;
\item [iii)] $w(x)=\sum_{n=1}^{\infty}{f_n(x)K_n(x)}$ converges absolutely and defines a continuous function of $x$.
\end{itemize}
Then $w(x)$ changes sign at most once and if the change occurs, it is from positive to negative.
\end{p1} 

\begin{p1}\label{le4}
For any positive integer $n$ and positive real number $x$, let $K_n(x)>0$ be a $TP_2$ function in $n\in \mathcal{N}^+$ and $x\in \mathcal{R}_+$. Assume that any function $f_n(x)$ be such that 
\begin{itemize}
\item [i)] for each $x\in \mathcal{R}_+$, $f_n(x)$ changes sign at most once and if the change occurs, it is from negative to positive as $n$ traverses in $\mathcal{N}^+$;
\item [ii)] for each $n\in \mathcal{N}^+$, $f_n(x)$ increases in $x$;
\item [iii)] $w(x)=\sum_{n=1}^{\infty}{f_n(x)K_n(x)}$ converges absolutely and defines a continuous function of $x$.
\end{itemize}
Then $w(x)$ changes sign at most once and if the change occurs, it is from negative to positive.
\end{p1} 
Now the question naturally arises,  whether any conclusion can be drawn about $w(x)$, for the combination of conditions other than those given in the Propositions \ref{le1}-\ref{le4}? Next one counterexample is provided to show that no conclusion regarding $w(x)$ can be drawn, as of the above mentioned propositions, for the combination of conditions given below:
\begin{itemize} 
\item[\textbf{Case-I}]   $K_n(x)$  is  RR2 in $n\in \mathcal{N}^+$ and $x \in \mathcal{R}_+$; $f_n(x)$ is increasing in  $x \in \mathcal{R}_+$ for all $n \in \mathcal{N}^+$;  and  $f_n(x)$ has at most one change of sign from negative to positive as $n$ traverses $\mathcal{N}^+.$
\item[\textbf{Case-II}]   $K_n(x)$  is  RR2 in $n\in \mathcal{N}^+$ and $x \in \mathcal{R}_+$;  $f_n(x)$ is decreasing in  $x \in \mathcal{R}_+$  for all $n \in \mathcal{N}^+$; and  $f_n(x)$ has at most one change of sign from  positive to negative as $n$ traverses $\mathcal{N}^+.$
\item[\textbf{Case-III}]   $K_n(x)$  is  TP2 function in $n\in \mathcal{N}^+$ and $x \in \mathcal{R}_+$;   $f_n(x)$ is decreasing in  $x \in \mathcal{R}_+$ for all $n \in \mathcal{N}^+$; and  $f_n(x)$ has at most one change of sign from  negative to positive as $n$ traverses $\mathcal{N}^+.$
\item[\textbf{Case-IV}]   $K_n(x)$  is  TP2 function in $n\in \mathcal{N}^+$ and $x \in \mathcal{R}_+$;  $f_n(x)$ is increasing in  $x \in \mathcal{R}_+$  for all $n \in \mathcal{N}^+$; and  $f_n(x)$ has at most one change of sign   positive to negative  as $n$ traverses $\mathcal{N}^+.$
\end{itemize}

 \begin{e1}
\textbf{Case-I:} Let us consider $K_n(x)=\frac{1}{nx^n}$, which is clearly RR2 function in $n\in \mathcal{N}^+$ and $x \in (1,\infty).$ Again consider, $f_n(x)=nx-5$, a function which is increasing in  $x \in (1,\infty).$ Again, for all $n\in \mathcal{N}^+$ it can be observed that it has at most one change of sign from  negative to positive as $n$ traverses $\mathcal{N}^+$, satisfying all conditions of Case I. Now it can be shown that 
\[
w_1(x)=\sum_{n=1}^{\infty}f_n(x)k_n(x)=\sum_{n=1}^{\infty}\frac{1}{x^{n-1}}-2\sum_{n=1}^{\infty}\frac{5}{nx^n}=\frac{x}{(x-1)}+5\ln\left(1-\frac{1}{x}\right)
\]
 is a continuous function on $x\in \mathcal{R}_+$, and from Figure $\ref{fig1.1}$ it is clear that $w_1(x)$  has more than one change of sign.\\

\begin{figure}[ht]
	\centering
	\includegraphics[height=5 cm]{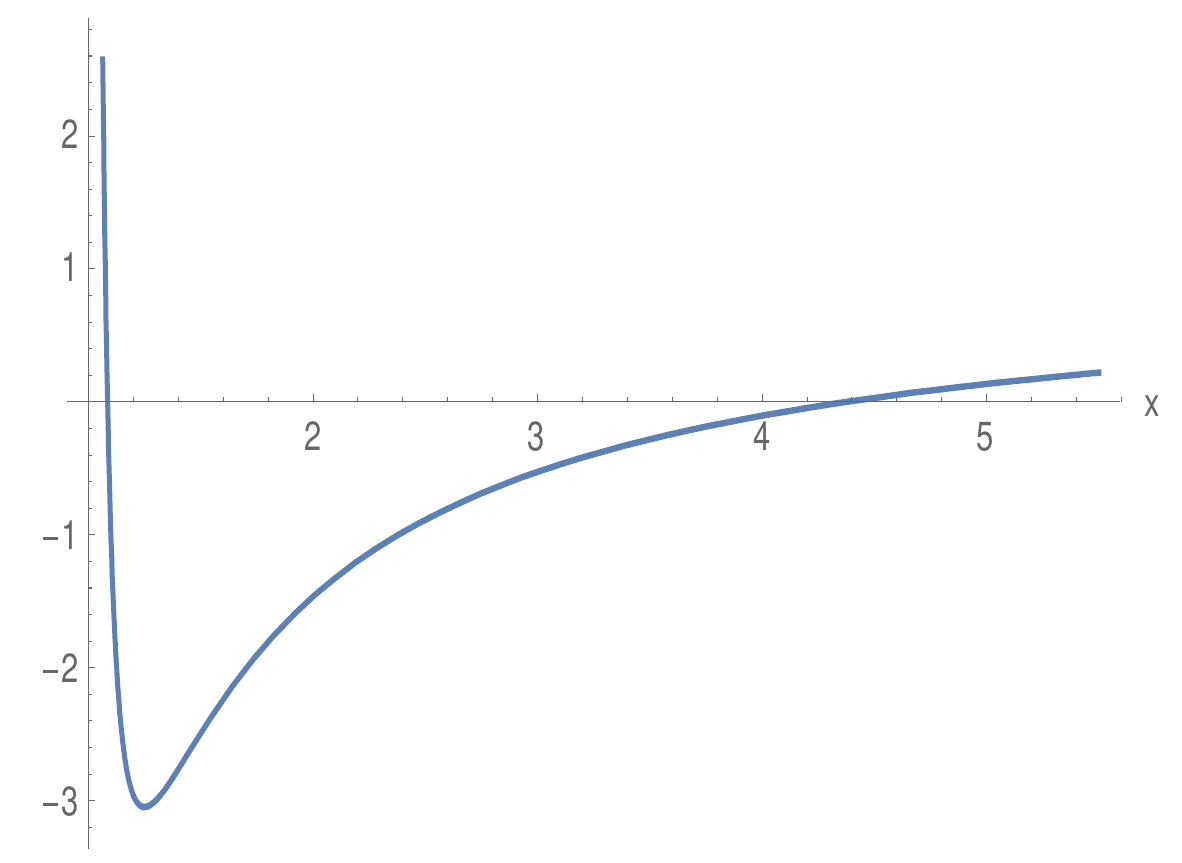}
	\caption{Graph of $w_1(x)$.} \label{fig1.1}
\end{figure}

\textbf{Case-II:} For the same $K_n(x)$ as of Case I, if $f_n(x)=5-nx$ is taken, which is a decreasing function in $x \in (1,\infty)$, then, for all $n\in \mathcal{N}^+$ it can be observed that it has at most one change of sign from  positive to negative as $n$ traverses $\mathcal{N}^+.$ Then clearly
\[
\sum_{n=1}^{\infty}f_n(x)k_n(x)=-w_1(x)
\]
 has more than one  change of sign which is evident from Figure $\ref{fig1.1}$.\\
\textbf{Case-III:}
Let $K_n(x) = \frac{x^n}{n}$  is a TP2 function in  in $n\in \mathcal{N}^+$ and $x \in (1,\infty).$ Also let  $f_n(x) = n - 10x$ which, clearly, for every  $n \in \mathcal{N}^+$ is decreeing in  $x \in (0,1)$. Additionally, for any positive real number $x$, $f_n(x)$ changes its sign at most once from negative to positive as $n$ varies across positive integers. Thus,

\[
w_2(x)=\sum_{n=1}^{\infty} f_n(x)k_n(x) = \sum_{n=1}^{\infty} {x^n} - 10x\sum_{n=1}^{\infty}\frac{x^n}{n}
\]
converges absolutely to the continuous function $\frac{x}{1-x}+10x\ln(1-x)$ for all $x \in (0,1).$
Now from Figure $\ref{fig1.2}$ it is clear that $w(x)$ has more than one change of sign, showing the result.
\begin{figure}[ht]
	\centering
	\begin{minipage}[b]{0.47\linewidth}
		\includegraphics[height=5 cm]{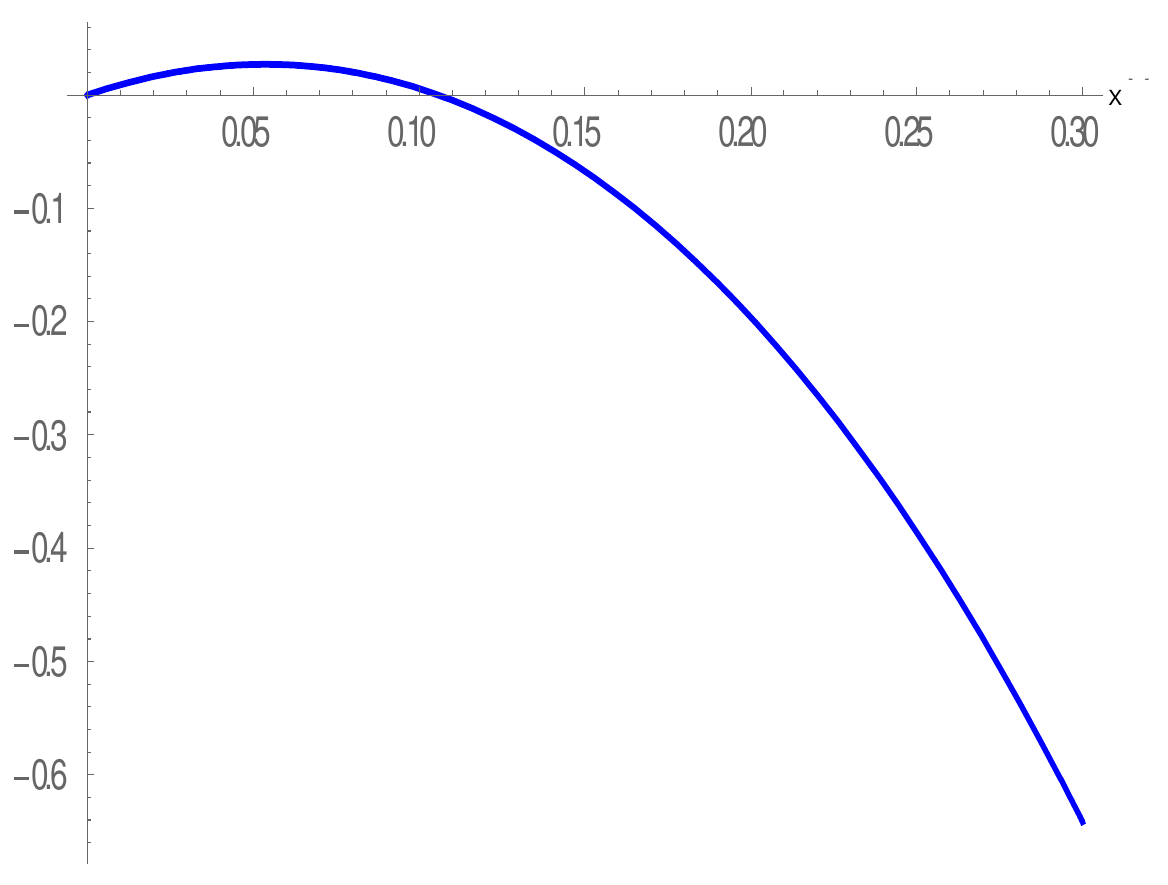}
		\centering{Graph of $w_2(x)$ on $[0,0.3].$  }
	\end{minipage}
	\quad
	\begin{minipage}[b]{0.48\linewidth}
		\includegraphics[height=5 cm]{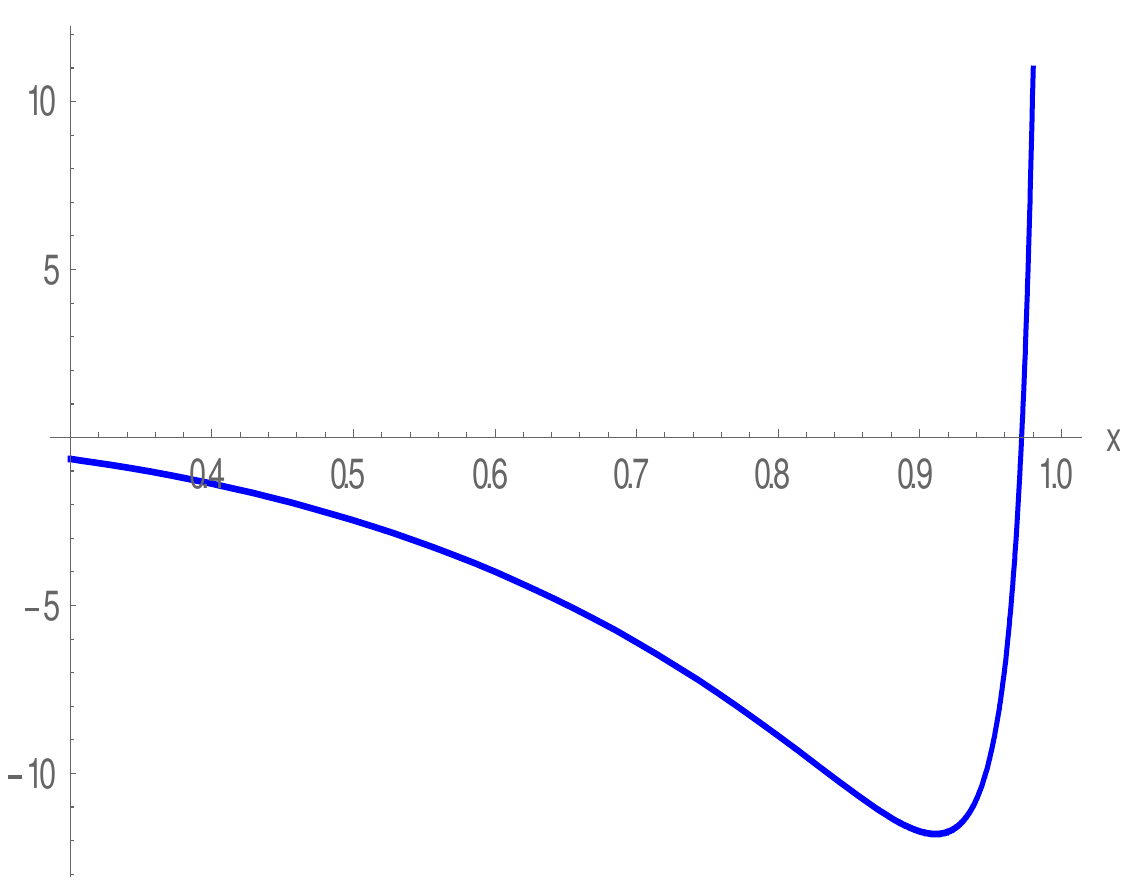}
		\centering{ Graph of $w_2(x)$ on $[0.3,1].$}
	\end{minipage}\caption{\label{fig1.2} Graph of $w_2(x)$.}
\end{figure}
\end{e1}

\textbf{Case-IV:}  For the same $K_n(x)$ as of Case III  if  $f_n(x)=10x-n$  is taken, which is increasing function in $x\in (0,1)$, then for all $n\in \mathcal{N}^+$ , it can be noticed that $f_n(x)$ changes its sign at most once from positive to negative as $n$ varies across positive integers. Then clearly
\[
\sum_{n=1}^{\infty} f_n(x)k_n(x) =-w_2(x) 
\]
is a continuous function in $x\in (0,1).$
So, from Figure $\ref{fig1.2}$, it is clear that $w_2(x)$  has more than one change of sign. 

\subsection{Comparisons Between Two Extreme Order Statistics}

 Let $X_i$ ($Y_i$), $i=1,2,\ldots,n$, be $n$ independent and non-identical random variables having distribution function $F_i\left(x\right)$ ($G_i\left(x\right)$), for all $i=1,2,\ldots,n$. Also let $F_{n:n}\left(x\right)$ and $G_{n:n}\left(x\right)$ ($\bar F_{1:n}\left(x\right)$ and $\bar G_{1:n}\left(x\right)$) be the distribution (survival) functions of $X_{n:n}$ and $Y_{n:n}$ ($X_{1:n}$ and $Y_{1:n}$) respectively with $X_{1:n}$ and $X_{n:n}$ being defined as $$ X_{1:n}=min\left\{X_1, X_2,...,X_n\right\}$$ and $$ X_{n:n}=max\left\{X_1, X_2,...,X_n\right\}.$$  Let $N$ is a discrete random variable with probability mass function (p.m.f) $p(n)$, $n\in\mathcal{N}^+$. Now, for the random variable $N$, if $F_{N:N}\left(x\right)$ and $G_{N:N}\left(x\right)$ ($\bar F_{1:N}\left(x\right)$ and $\bar G_{1:N}\left(x\right)$) denote the distribution (survival) functions of $X_{N:N}$ and $Y_{N:N}$ ($X_{1:N}$ and $Y_{1:N}$) respectively, then it can be written that 
 \begin{equation}\label{eq2}
 	F_{N:N}(x)=\sum_{n=1}^\infty F_{n:n}(x)p(n),\ G_{N:N}(x)=\sum_{n=1}^\infty G_{n:n}(x)p(n),
 \end{equation}
 \begin{equation}\label{eq3}
 	\bar F_{1:N}(x)=\sum_{n=1}^\infty \bar F_{1:n}(x)p(n)\ \text{and}\ \bar G_{1:N}(x)=\sum_{n=1}^\infty \bar G_{1:n}(x)p(n).
 \end{equation}

The following theorems show that under certain restrictions, if there exists hazard rate (reversed hazard rate) ordering between $X_{1:n}$ and $Y_{1:n}$ ($X_{n:n}$ and $Y_{n:n}$), then the same ordering also exists between $X_{1:N}$ and $Y_{1:N}$ ($X_{N:N}$ and $Y_{N:N}$).

\begin{t1}\label{th1}
Suppose $X_{1:n}\sim \bar F_{1:n}(x)$ and  $Y_{1:n}\sim \bar G_{1:n}(x).$ Also let $N$ be a random variable having pmf $p(n).$ Now, for all $x\geq 0$, if $\frac{\bar{F}_{1:n}(x)}{\bar{G}_{1:n}(x)}$ is increasing in $n\in\mathcal{N}^+$, then $X_{1:n}\leq_{hr} Y_{1:n}$ implies $X_{1:N}\leq_{hr} Y_{1:N}.$
	\end{t1}
	{\bf Proof:} Here, we need to prove that $\frac{\bar{F}_{1:N}(x)}{\bar{G}_{1:N}(x)}$ is decreasing in $x$ or equivalently $\frac{\sum_{n=1}^{\infty}{\bar{F}_{1:n}(x)p(n)}}{\sum_{n=1}^{\infty}{\bar{G}_{1:n}(x)p(n)}}$ is decreasing in $x.$ Now for any $\lambda>0$, let 
	\begin{eqnarray*}
	w_3(x) &=&\sum_{n=1}^{\infty}{\left[\bar{F}_{1:n}(x)-\lambda\bar{G}_{1:n}(x)\right]p(n)}\\
	     &=& \sum_{n=1}^{\infty}{\bar{G}_{1:n}(x)\left[\frac{\bar{F}_{1:n}(x)}{\bar{G}_{1:n}(x)}-\lambda\right]p(n)}\\
			&=& \sum_{n=1}^{\infty}{\bar{G}_{1:n}(x) f_n(x) p(n)} \ \text{(say),}
	\end{eqnarray*}
where $f_n(x)=\left[\frac{\bar{F}_{1:n}(x)}{\bar{G}_{1:n}(x)}-\lambda\right]$.\\
Now, $X_{1:n}\leq_{hr} Y_{1:n}$ implies that $\frac{\bar{F}_{1:n}(x)}{\bar{G}_{1:n}(x)}$ is decreasing in $x$, and hence $f_n(x)$ is also decreasing in $x.$ \\
Since $\frac{\bar{F}_{1:n}(x)}{\bar{G}_{1:n}(x)}$ is increasing in $n$, it is clear that $f_n(x)$ has at most one change of sign as $n$ traverses and if it does occur, it is from negative to positive.\\ 
Also, Theorem \ref{theo11} concludes that for any two positive integers $n_1\leq n_2$, $Y_{1:n_1}\geq_{hr}Y_{1:n_2}$, which by definition of $hr$ ordering indicates that $\frac{\bar{G}_{1:n_1}(x)}{\bar{G}_{1:n_2}(x)}$ is increasing in $x$, yielding that $\bar{G}_{1:n}(x)$ is $RR_2$ function in $n$ and $x$.\\
Thus, by Proposition $\ref{le1}$, $w_3(x)$ changes it's sign at most once and if the changes occurs, it is from positive to negative. Now, as we know that (see proof of Lemma 3.2 of Nanda and Das~\cite{nan}), for any $x\geq 0$, $f(x)-cg(x)$ has at most one change of sign from  positive to negative can equivalently  written as $\frac{f(x)}{g(x)}$ is decreasing in $x\in (0,\infty)$, thus it can be concluded that $\frac{\sum_{n=1}^{\infty}{\bar{F}_{1:n}(x)}p(n)}{\sum_{n=1}^{\infty}\bar{G}_{1:n}(x)p(n)}$ is decreasing in $x$, proving the result.$\hfill\Box$\\
Next, we illustrate Theorem \ref{th1} with two examples. In the first example, independent heterogenous location family distributed random variable is assumed while $inid$ Weibull distributed random variable is chosen in the second example.\\ 

{\bf Example 1}:
Let, $X\sim \bar F(x)$ is a random variable having decreasing failure rate function $r(x)$. Also let, for $n\in\mathcal{N}^+$, $X_i\sim \bar F\left(x-\lambda_i\right)$ and $Y_i\sim \bar F\left(x-\mu_i\right),~i=1,2,...,n$ are two sequences of $inid$ random variables where $\lambda_i\geq\mu_i$. Then, the survival and failure rate functions of $X_{1:n}$ and $Y_{1:n}$ are obtained as  
\begin{equation}\label{eq4}
\bar F_{1:n}(x)=\prod_{i=1}^n\bar F\left(x-\lambda_i\right),\ \bar G_{1:n}(x)=\prod_{i=1}^n\bar F\left(x-\mu_i\right),
\end{equation}
and 
\begin{equation}\label{eq5}
r_{1:n}(x)=\sum_{i=1}^n r\left(x-\lambda_i\right)\ \text{and}\ s_{1:n}(x)=\sum_{i=1}^n r\left(x-\mu_i\right)
\end{equation}
respectively. \\
Since, $\lambda_i\geq\mu_i$ and $r(x)$ is decreasing in $x$, it can be written that $r\left(x-\lambda_i\right)\geq r\left(x-\mu_i\right)$ for $i=1,2,\ldots,n.$ Hence, from (\ref{eq5}), it is clear that
$$r_{1:n}(x)=\sum_{i=1}^n r\left(x-\lambda_i\right)\geq \sum_{i=1}^n r\left(x-\mu_i\right)=s_{1:n}(x),$$
proving that $X_{1:n}\leq_{hr}Y_{1:n}$.\\
Now, for all $n\in\mathcal{N}^+,$ $\lambda_i\geq\mu_i$ gives $\bar{F}\left(x-\lambda_i\right)\geq \bar{F}\left(x-\mu_i\right).$ From (\ref{eq4}), it can be immediately concluded that   
$$\frac{\bar F_{1:n}(x)}{\bar G_{1:n}(x)}=\prod_{i=1}^n \frac{\bar{F}\left(x-\lambda_i\right)}{\bar{F}\left(x-\mu_i\right)}$$
is increasing in $n$.\\
Furthermore, from (\ref{eq4}) it can be written that for all $n_1\leq n_2$, 
\begin{eqnarray*}
\frac{\bar G_{1:n_1}(x)}{\bar G_{1:n_2}(x)}=\frac{\prod_{i=1}^{n_1}\bar F\left(x-\mu_i\right)}{\prod_{i=1}^{n_2}\bar F\left(x-\mu_i\right)}=\frac{1}{\prod_{i=n_1+1}^{n_2}\bar F\left(x-\mu_i\right)},
\end{eqnarray*}
which is clearly increasing in $x$. Thus all the conditions of Theorem \ref{th1} are satisfied. Hence, it can be concluded that $X_{1:N}\leq_{hr}Y_{1:N}$.\\

\begin{figure}[ht]
\centering
\includegraphics[height=7 cm]{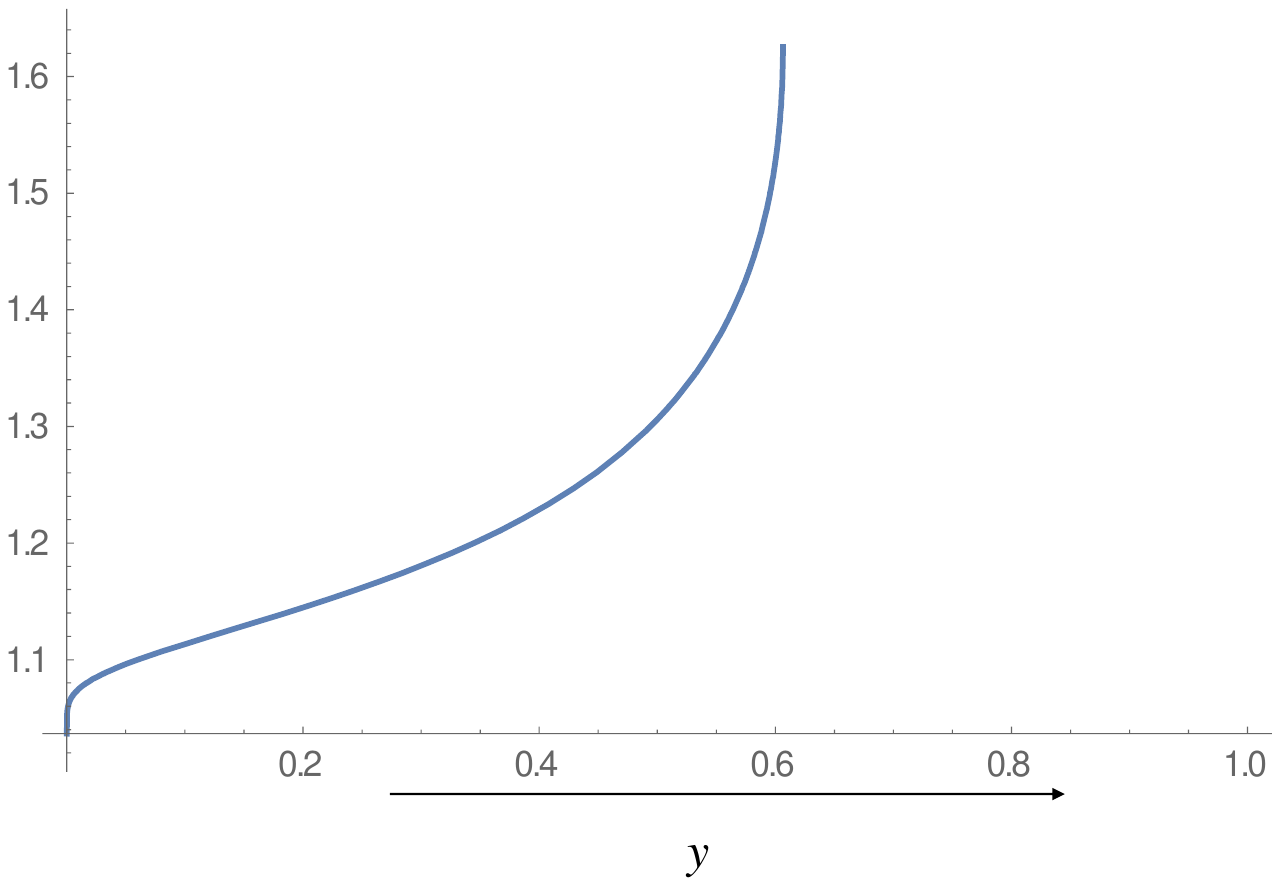}
\caption{Graph of $\frac{\bar{F}_{1:N}(-\ln y)}{\bar{G}_{1:N}(-\ln y)}$} \label{figure1}
\end{figure}

{\bf Example 5}: Let us assume that $X$ follows Weibull distribution with survival function $$\bar F(x)=e^{-2x^{0.5}}.$$ Clearly the failure rate function of $X$, $r(x)=\frac{1}{\sqrt x}$ is decreasing in $x.$\\ For $i=1,2\ldots, 5$, suppose $X_i\sim \bar F\left(x-\lambda_i\right)$ and $Y_i\sim \bar F\left(x-\mu_i\right)$ where\\ $\left(\lambda_1,\lambda_2,\lambda_3,\lambda_4,\lambda_5\right)=(0.1,0.2,0.3,0.4,0.5)$ and $\left(\mu_1,\mu_2,\mu_3,\mu_4,\mu_5\right)=(0.05,0.15,0.25,0.35,0.45)$, which implies that $\lambda_i\geq\mu_i$. So, for $n=1,2\ldots, 5$, it can be written that\\ $\bar{F}_{1:n}(x)=\prod_{i=1}^n \bar F(x-\lambda_i)$ and  $\bar{G}_{1:n}(x)=\prod_{i=1}^n \bar F(x-\mu_i)$. Now, if it is assumed that  $$X_1, X_2, X_3~(Y_1, Y_2, Y_3)~\text{are selected with probability} P(N=3)=p(3)=1/5,$$ $$X_1, X_2, X_3, X_4 ~(Y_1, Y_2, Y_3,Y_4)~\text{are selected with probability} P(N=4)=p(4)=2/5~\text{and}$$ $$X_1, X_2, X_3, X_4, X_5~(Y_1, Y_2, Y_3, Y_4, Y_5)~\text{are selected with probability} P(N=5)=p(5)=2/5,$$ then from Figure \ref{figure1}, it can be verified that
$$\frac{\bar{F}_{1:N}(x)}{\bar{G}_{1:N}(x)}=\frac{{\bar{F}_{1:3}(x)}p(3)+{\bar{F}_{1:4}(x)}p(4)+{\bar{F}_{1:5}(x)}p(5)}{{\bar{G}_{1:3}(x)}p(3)+{\bar{G}_{1:4}(x)}p(4)+{\bar{G}_{1:5}(x)}p(5)}$$ is decreasing in $x$, proving the result. Here, substitution $x=-\ln y$, $0\leq y\leq 1$ is used to capture the real line.
 \begin{r1}
It is to be noted that Theorem \ref{th1} is true for usual stochastic ($st$) ordering under no restriction. If there exists usual stochastic ($st$) ordering between $X_{1:n}$ and $Y_{1:n},$ there will also exist usual $st$ ordering between $X_{1:N}$ and $Y_{1:N}$. The proof is trivial. The same can be observed in Theorems \ref{th2}-\ref{th4} under no additional sufficient conditions.
\end{r1}
The next few theorems can be proved in the same fashion as above and hence the details are omitted.
\begin{t1}\label{th2}
Suppose $X_{1:n}\sim \bar F_{1:n}(x)$ and  $Y_{1:n}\sim \bar G_{1:n}(x).$ Also let $N$ be a random variable having pmf $p(n).$ Now, for all $x\geq 0$, if $\frac{\bar{F}_{1:n}(x)}{\bar{G}_{1:n}(x)}$ is decreasing in $n\in\mathcal{N}^+$, then $X_{1:n}\geq_{hr} Y_{1:n}$ implies $X_{1:N}\geq_{hr} Y_{1:N}.$
	\end{t1}	
{\bf Proof}: The theorem can be proved using Proposition \ref{le2}, Theorem \ref{theo11} and following the proof of Theorem \ref{th1}.
	
	\begin{t1}\label{th3}
	Suppose $X_{n:n}\sim \bar F_{n:n}(x)$ and  $Y_{n:n}\sim \bar G_{n:n}(x).$ Also let $N$ be a random variable having pmf $p(n).$ Now, for all $x\geq 0$, if $\overline r_{Y_{n:n}}(x)$, the reversed hazard rate function of $Y_{n:n}$ is increasing in $n\in\mathcal{N}^+$ and $\frac{F_{n:n}(x)}{G_{n:n}(x)}$ is decreasing in $n\in\mathcal{N}^+$, then $X_{n:n}\leq_{rh} Y_{n:n}$ implies $X_{N:N}\leq_{rh} Y_{N:N}.$
	\end{t1}	
	{\bf Proof}: If $\overline r_{Y_{n:n}}(x)$ is increasing in $n$, then for any two positive integers $n_1\leq n_2$ and for all $x\geq 0$ it can be written that $\overline r_{Y_{n_1:n_1}}(x)\leq \overline r_{Y_{n_2:n_2}}(x)$, which can equivalently conclude that $\frac{G_{n_1:n_1}(x)}{G_{n_2:n_2}(x)}$ is decreasing in $x\geq 0$. Thus, the theorem can be proved using Proposition~$\ref{le3}$ and following the line of the proof of Theorem \ref{th1}.
	
	\begin{t1}\label{th4}
	Suppose $X_{n:n}\sim \bar F_{n:n}(x)$ and  $Y_{n:n}\sim \bar G_{n:n}(x).$ Also let $N$ be a random variable having pmf $p(n).$ Now, for all $x\geq 0$, if $\overline r_{Y_{n:n}}(x)$, the reversed hazard rate function of $Y_{n:n}$ is increasing in $n\in\mathcal{N}^+$ and $\frac{F_{n:n}(x)}{G_{n:n}(x)}$ is increasing in $n\in\mathcal{N}^+$, then $X_{n:n}\geq_{rh} Y_{n:n}$ implies $X_{N:N}\geq_{rh} Y_{N:N}.$
	\end{t1}
{\bf Proof}: The theorem can be proved using Proposition~\ref{le4} and following the proof of Theorem \ref{th3}. 
\begin{r1}
It is to be noted that Theorem \ref{th3} and \ref{th4} are proved with an additional condition that ``for all $x\geq 0$, the reversed hazard rate function ($\overline r_{Y_{n:n}}(x)$) of $Y_{n:n}$ is increasing in $n$", while no such condition is required to prove Theorems \ref{th1} and \ref{th2}. This can be rationalized by  Theorem 1.B.28 of~\cite{shak} (see Theorem \ref{theo11} above) which states that for all $x\geq 0$, the hazard rate function ($r_{Y_{1:n}}(x)$) of $Y_{1:n}$ is always increasing in $n$ which is not always true for reversed hazard rate function. In fact, Theorem 1.B.57 of~\cite{shak} (see Theorem \ref{theo12} above) uses more restrictive sufficient conditions $X_m\leq_{rh}X_j$ for all $j=1,2,\ldots, m-1$ to prove similar result. Hence, one additional condition is used in Theorems \ref{th3}~and~\ref{th4} to avoid such restrictive sufficient conditions. 
\end{r1}

The subsequent theorems demonstrate that when specific conditions are met, if there is a likelihood-ratio ($lr$) dominance between $X_{1:n}$ and $Y_{1:n}$ ($X_{n:n}$ and $Y_{n:n}$). Then the dominance is preserved for $X_{1:N}$ and $Y_{1:N}$ $(X_{N:N}\; \text{and} \; Y_{N:N})$.

\begin{t1}\label{th5}
	
	Suppose $X_{1:n}$ and  $Y_{1:n}$  be two minimum order statistic having density functions $f_{1:n}(x)$ and $g_{1:n}(x)$ respectively. Also let $N$ be a discrete random variable having p.m.f $p(n).$ Now, for all $x\geq 0$, if $\frac{g_{1:n}(x)}{f_{1:n}(x)}$, is increasing in $n$ and $n_1\leq n_2$ implies $X_{1:n_1}\geq_{lr} X_{1:n_2}$. Then $X_{1:n}\geq_{lr} Y_{1:n}$ implies $X_{1:N}\geq_{lr} Y_{1:N}.$
\end{t1}

{\bf Proof:} Here, we need to prove that $\frac{g_{1:N}(x)}{f_{1:N}(x)}$ is decreasing in $x$ or equivalently $\frac{\sum_{n=1}^{\infty}{g_{1:n}(x)p(n)}}{\sum_{n=1}^{\infty}{f_{1:n}(x)p(n)}}$ is decreasing in $x.$ Now for any $\lambda>0$, let 
\begin{eqnarray*}
w_4(x) &=&\sum_{n=1}^{\infty}{\left[g_{1:n}(x)-\lambda f_{1:n}(x)\right]p(n)}\\
&=& \sum_{n=1}^{\infty}{f_{1:n}(x)\left[\frac{g_{1:n}(x)}{f_{1:n}(x)}-\lambda\right]p(n)}\\
&=& \sum_{n=1}^{\infty}{f_{1:n}(x) f_n(x) p(n)} \ \text{(say),}
\end{eqnarray*}
where $f_n(x)=\left[\frac{g_{1:n}(x)}{f_{1:n}(x)}-\lambda\right]$.\\
Now, $X_{1:n}\geq_{lr} Y_{1:n}$ implies that $\frac{g_{1:n}(x)}{f_{1:n}(x)}$ is decreasing in $x$, and hence $f_n(x)$ is also decreasing in $x.$ \\
Since $\frac{g_{1:n}(x)}{f_{1:n}(x)}$ is increasing in $n$, it is clear that $f_n(x)$ has at most one change of sign as $n$ traverses and if it does occur, it is from negative to positive.\\
Also, as for any two positive integers $n_1\leq n_2$, $X_{1:n_1}\geq_{lr}X_{1:n_2}$, then $\frac{f_{1:n_1}(x)}{f_{1:n_2}(x)}$ is increasing in $x$, yielding that $f_{1:n}(x)$ is $RR_2$ function in $n$ and $x$.\\ 
Thus, by Proposition $\ref{le1}$, $w_4(x)$ changes it's sign at most once and if the changes occurs, it is from positive to negative, proving the result.$\hfill\Box$\\

Next one example is provided to illustrate Theorem~\ref{th5}.\\

{\bf Example 6}: Let us assume that \( X \) follows a Weibull distribution with survival function:
	\[ \bar F(x) = e^{-0.2x^{0.6}}. \]
	
	For \( i=1,2,\ldots,5 \), suppose \( X_i\sim \bar F(x-\alpha_i) \) and \( Y_i\sim \bar F(x-\beta_i) \) where
	\[(\alpha_1,\alpha_2,\alpha_3,\alpha_4,\alpha_5) = (0.05,0.15,0.25,0.35,0.45)\]
	and
	\[ (\beta_1,\beta_2,\beta_3,\beta_4,\beta_5) = (0.1,0.2,0.3,0.4,0.5). \] So, for \( n=1,2,\ldots, 5 \), it can be written that
	$ \bar{F}_{1:n}(x) = \prod_{i=1}^n \bar F(x-\alpha_i) $ and
	$ \bar{G}_{1:n}(x) = \prod_{i=1}^n \bar F(x-\beta_i). $
	
\begin{figure}[ht]
		\centering
		\includegraphics[height=5 cm]{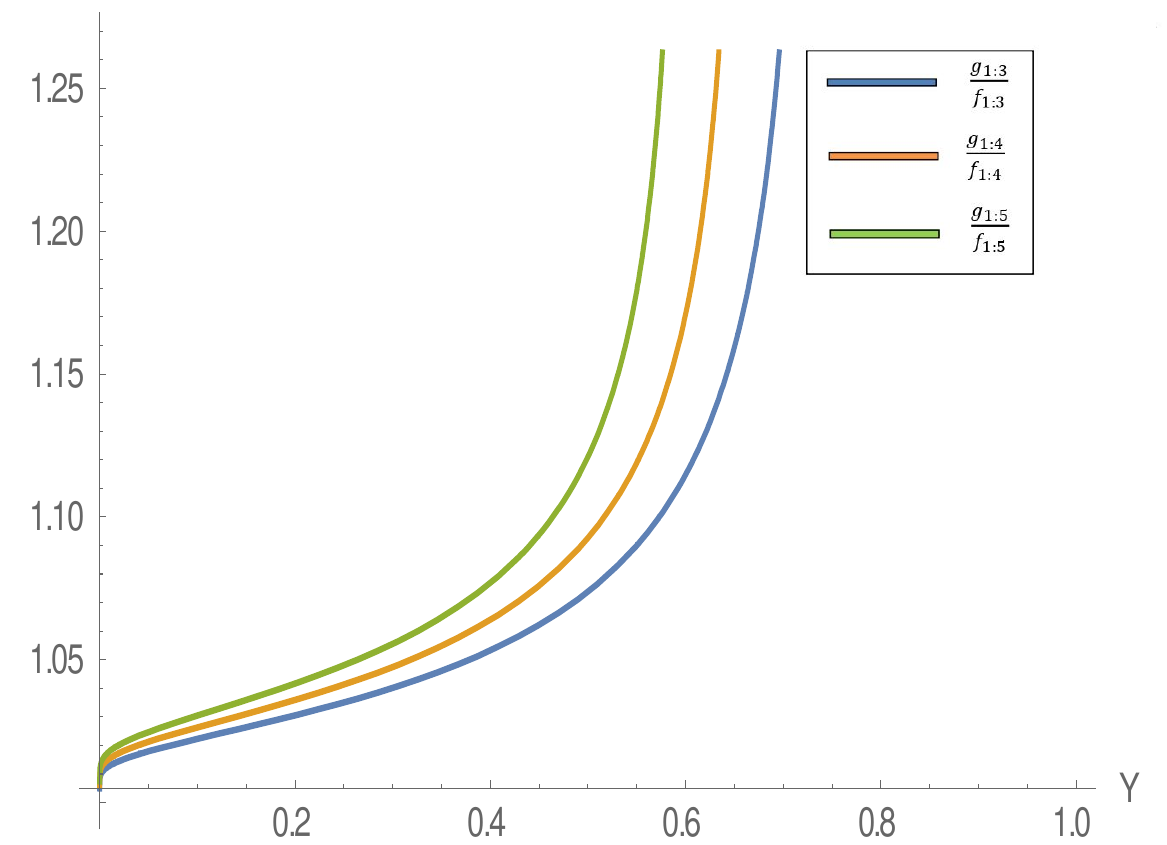}
		\caption{Graph of \( \frac{g_{1:n}(y)}{f_{1:n}(y)} \)} \label{figure2.1}
	\end{figure}

\begin{figure}[ht]
		\centering
		\includegraphics[height=5 cm]{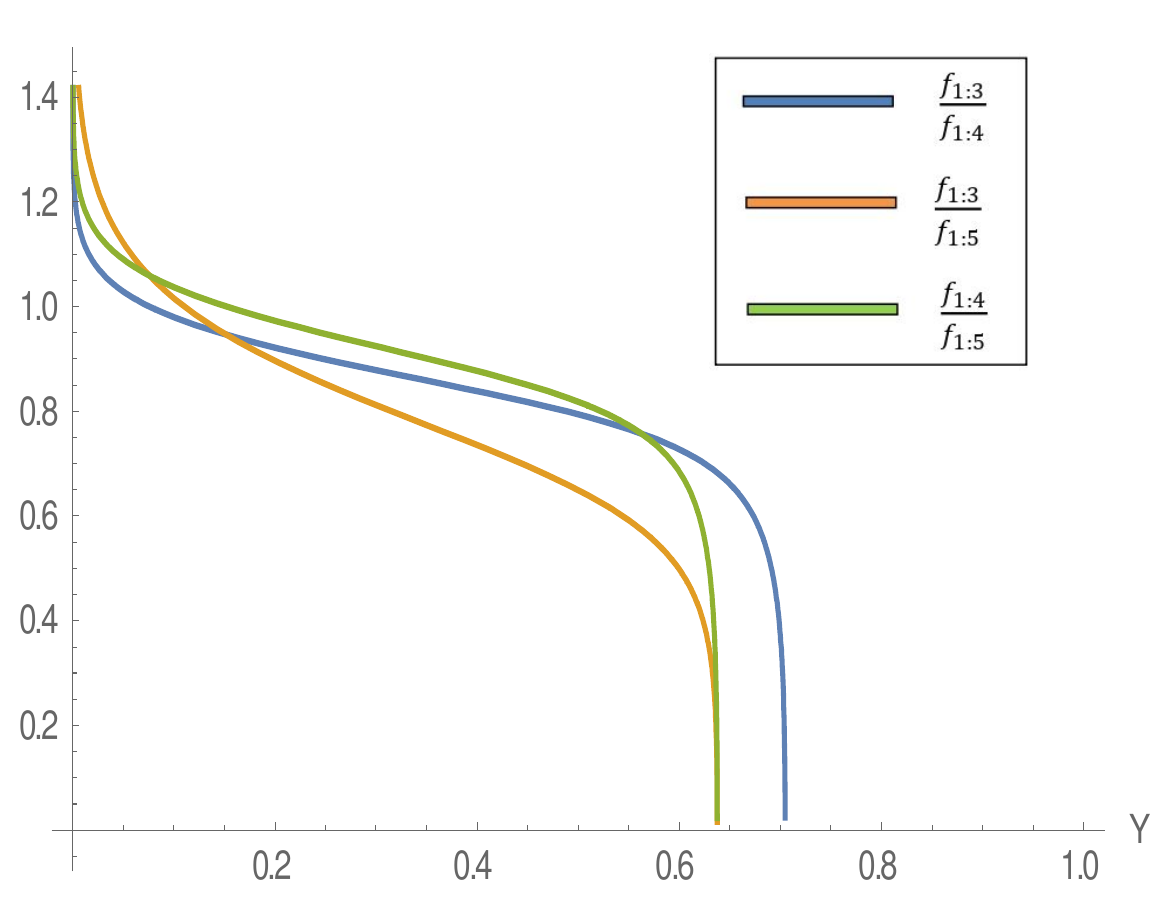}
		\caption{Graph of  $\frac{f_{1:n_1}(y)}{f_{1:n_2}(y).}$} \label{figure2.2}
	\end{figure}

\begin{figure}[ht]
		\centering
		\includegraphics[height=5 cm]{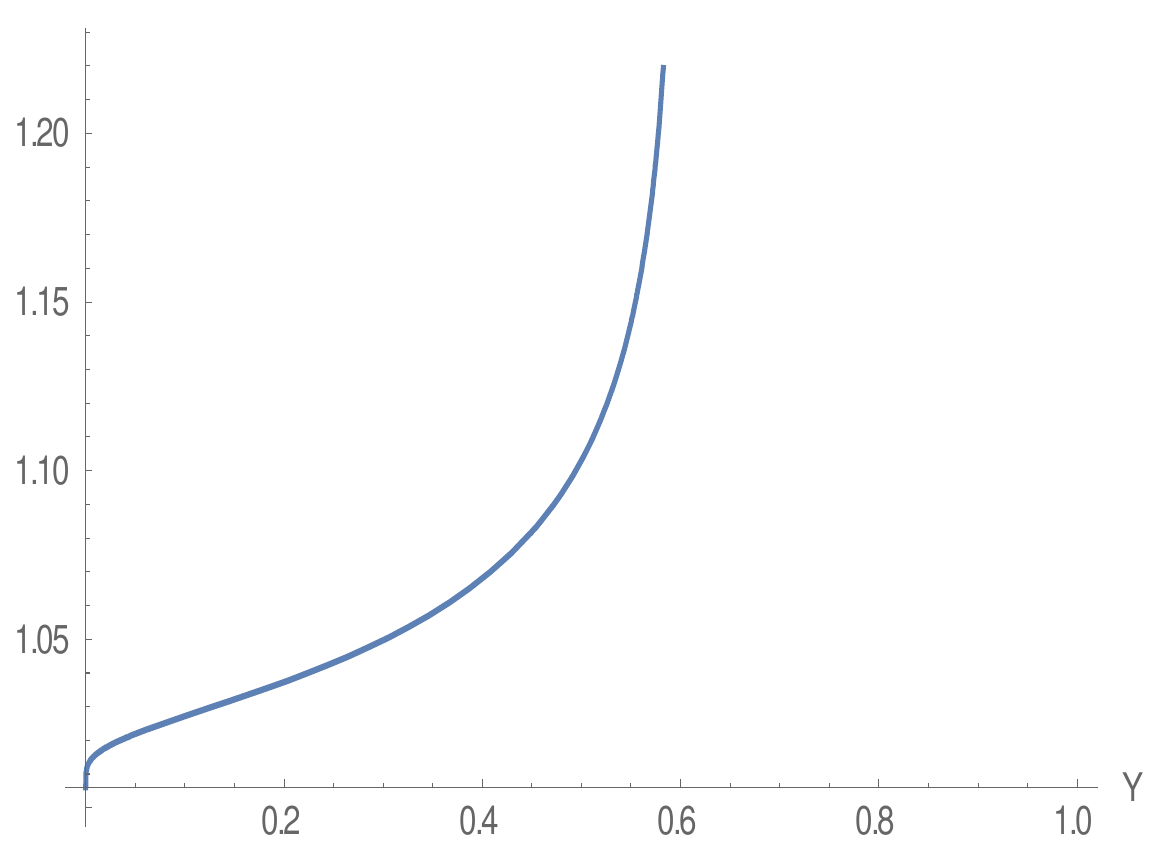}
		\caption{Graph of \( \frac{g_{1:N}(y)}{f_{1:N}(y)} \)} \label{figure2.3}
	\end{figure}

	Now from Figure~\ref{figure2.1}, it is clear that for all \( x\geq 0 \), if \( \frac{g_{1:n}(x)}{f_{1:n}(x)} \) is decreasing in $x$ and increasing in $n$ . Again, from Figure~\ref{figure2.2} it is also clear that for all $3\leq n_1 \leq n_2 \leq 5,$ $X_{1:n_1} \geq_{lr} .
	X_{1:n_2}.$
	Now, if it is assumed that
	\[ X_1, X_2, X_3~(Y_1, Y_2, Y_3)~\text{are selected with probability} P(N=3)=p(3)=\frac{1}{5}, \]
	\[ X_1, X_2, X_3, X_4 ~(Y_1, Y_2, Y_3,Y_4)~\text{are selected with probability} P(N=4)=p(4)=\frac{2}{5}, \]
	and
	\[ X_1, X_2, X_3, X_4, X_5~(Y_1, Y_2, Y_3, Y_4, Y_5)~\text{are selected with probability} P(N=5)=p(5)=\frac{2}{5}, \]
	then from Figure~\ref{figure2.3}, it can be verified that
	\[
	\frac{g_{1:N}(x)}{f_{1:N}(x)} = \frac{{g_{1:3}(x)}p(3)+{g_{1:4}(x)}p(4)+{g_{1:5}(x)}p(5)}{{f_{1:3}(x)}p(3)+{f_{1:4}(x)}p(4)+{f_{1:5}(x)}p(5)}
	\]
	is decreasing in \( x \), proving the result. Here, substitution \( y=-\ln x \), \( 0\leq y\leq 1 \) is used for each of Figure~\ref{figure2.1}-\ref{figure2.3}, to capture the whole real line.
	
\begin{t1}\label{th6}
	
	Suppose $X_{1:n}$ and  $Y_{1:n}$  be two minimum order statistic having density functions $f_{1:n}(x)$ and $g_{1:n}(x)$ respectively. Also let $N$ be a  discrete random variable having pmf $p(n).$ Now, for all $x\geq 0$, if $\frac{g_{1:n}(x)}{f_{1:n}(x)}$ is decreasing in $n$, and $n_1\leq n_2$ implies $X_{1:n_1}\geq_{lr} X_{1:n_2}$. Then $X_{1:n}\leq_{lr} Y_{1:n}$ implies $X_{1:N}\leq_{lr} Y_{1:N}.$
\end{t1}
{\bf Proof}: The theorem can be proved using Proposition \ref{le2},  and following the proof of Theorem \ref{th5}.
\begin{t1}\label{th7}
	
	Suppose $X_{n:n}$ and  $Y_{n:n}$  be two maximum order statistic having density functions $f_{n:n}(x)$ and $g_{n:n}(x)$ respectively. Also let $N$ be a discrete random variable having p.m.f $p(n).$ Now, for all $x\geq 0$, if $\frac{g_{n:n}(x)}{f_{n:n}(x)}$ is increasing in $n$ and $n_1\leq n_2$ implies $X_{n_1:n_1}\leq_{lr} X_{n_2:n_2}$. Then $X_{n:n}\leq_{lr} Y_{n:n}$ implies $X_{N:N}\leq_{lr} Y_{N:N}.$
\end{t1}
{\bf Proof}: Since for any two positive integer $n_1\leq n_2$ implies $X_{n_1:n_1}\leq_{lr} X_{n_2:n_2}$, which can be equivalently conclude that $\frac{f_{n_1:n_1}}{f_{n_2:n_2}}$ is decreasing in $x$ i.e. $f_{n:n}$ is $TP2$ in $n$ and $x$. Thus, the theorem can be proved using Proposition $\ref{le4}$  and  following the line of the proof of Theorem $\ref{th5}$. 
 \begin{t1}\label{th8}
	
Suppose $X_{n:n}$ and  $Y_{n:n}$  be two maximum order statistic having density functions $f_{n:n}(x)$ and $g_{n:n}(x)$ respectively. Also let $N$ be a discrete random variable having p.m.f $p(n).$ Now, for all $x\geq 0$, if $\frac{g_{n:n}(x)}{f_{n:n}(x)}$ is decreasing in $n$ and $n_1\le n_2$ implies $X_{n_1:n_1}\leq_{lr} Y_{n_2:n_2}$. Then $X_{n:n}\geq_{lr} X_{n:n}$ implies $X_{N:N}\geq_{lr} Y_{N:N}.$
\end{t1}
{\bf Proof}: The theorem can be proved using Proposition \ref{le3} and following the proof of Theorem $\ref{th7}$.

\section{Conclusion}
\setcounter{equation}{0}   
	\hspace*{0.3 in} The paper aims to provide some new results on the comparison of extremes with random number of independent and non-identically distributed random variables. The variation diminishing property (~\cite{ka}) is extended to all possible restrictions and the associated propositions are proved which are used to derive the key results. Stochastic comparison of $k$-out-of-$n$ systems where the components are non-identically distributed and the number of components is random, can be an interesting problem for future research.  
\vskip2pt

\noindent{\bf Statements and Declarations}\\
No potential conflict of interest was reported by the author(s).

\section*{Acknowledgement}
A fruitful discussion with Asok K Nanda, Professor of Indian Institute of Science and Educational Research, Kolkata, India is gratefully acknowledged. We also thank the anonymous reviewers and the Editor for their constructive comments which led to the improved version of the manuscript.

\end{document}